\documentstyle[12pt]{article}

\newtheorem{lemma}{{\sc Lemma}}[section]
\newtheorem{theo}{{\sc Theorem}}

\newtheorem{defn}{{\em Definition}}

\newcommand{\proof}{{\sc Proof\,:}}

\begin{document}


\title{Contractible Hamiltonian Cycles in Triangulated Surfaces}

\author{
{Ashish Kumar Upadhyay}\\[2mm]
{\normalsize Department of Mathematics}\\ {Indian Institute of Technology Patna }  \\
{\normalsize Patliputra Colony, Patna 800\,013,  India.}\\[1mm]
{\small upadhyay@iitp.ac.in} }
\maketitle

\vspace{-5mm}

\hrule

\begin{abstract} A triangulation of a surface is called $q$-equivelar if each of its vertices
is incident with exactly $q$ triangles. In 1972 Altshuler had shown that an equivelar triangulation of torus
has a Hamiltonian Circuit. Here we present a necessary and sufficient condition for
existence of a contractible Hamiltonian Cycle in equivelar triangulation of a
surface.
\end{abstract}

{\small

{\bf AMS classification\,:} 57Q15, 57M20, 57N05.

{\bf Keywords\,:} Contractible Hamiltonian cycles, Proper Trees in Maps,
Equivelar \newline \mbox{} \hspace{25mm} Triangulations.
}

\bigskip

\hrule

\section{Introduction}

A graph $G := (V, E)$ is without loops and such that no more than one edge joins two vertices.
A {\em map} on a surface $S$ is an embedding of a graph $G$ with finite
number of vertices such that the components of $S \setminus G$ are topological $2$-cells.
Thus, the closure of a cell in $S\setminus G$ is a $p-gonal$ disk, $i.e.$ a 2-disk whose
boundary is a $p-$gon for some integer $p \geq 3$.

A map is called {\em $\{p, q\}$ equivelar} if each vertex is incident with exactly $q$
numbers of $p$-gons. If $p = 3$ then the map is called a $q$ - equivelar triangulation or a
degree - regular triangulation of type $q$. A map is called a {\em Simplicial Complex} if
each of its faces is a simplex. Thus a triangulation is a {\em Simplicial Complex}
For a simplicial complex K, the graph consisting of the edges and vertices of
$K$ is called the {\em edge-graph} of $K$ and is denoted by $EG(K)$.

If $X$ and $Y$ are two simplicial complexes, then a (simplicial) isomorphism from $X$ to
$Y$ is a bijection $\phi : V (X)\rightarrow V (Y)$ such that for $\sigma \subset V (X)$,
$\sigma$ is a simplex of $X$ if and only if $\phi(\sigma)$ is a simplex of $Y$.
Two simplicial complexes $X$ and  $Y$ are called (simplicially)
isomorphic (and is denoted by $X\cong Y$ ) when such an isomorphism exists. We identify
two complexes if they are isomorphic. An isomorphism from a simplicial complex $X$ to
itself is called an automorphism of $X$. All the automorphisms of $X$ form a group, which
is denoted by $Aut(X)$.


In 1956, Tutte \cite{Tutte} showed that every 4-connected planar graph has a Hamiltonian
cycle. Later in 1970, Gr$\ddot{u}$nbaum conjectured that every 4-connected graph
which admits an embedding in the torus has a Hamiltonian cycle. In the same article
he also remarked that - probably there is a function $c(k)$ such that each $c(k)$-connected
graph of genus at most $k$ is Hamiltonian.

In 1972, Duke \cite{duke} showed the existence of such a function and gave an estimate
$[\frac{1}{2}(5 + \sqrt{16\,k + 1})]
\leq c(k) \leq \{3 + \sqrt{6\,k + 3}\}$ where $k\geq 1$.

A. Altshuler \cite{a1}, \cite{a2} studied Hamiltonian cycles and paths in the edge
graphs of equivelar maps on the torus. That is in the maps which are equivelar of
types $\{3, 6\}$ and $\{4, 4\}$. He showed that in the graph consisting of
vertices and edges of equivelar maps of above type there exists a Hamiltonian cycle.
He also showed that a Hamiltonian cycle exists in every 6-connected graph on the torus.

In 1998, Barnette \cite{barnette3} showed that any 3-connected graph other than $K_4$
or $K_5$ contains a contractible cycle or contains a simple configuration as
subgraphs.

In this article we present a necessary and
sufficient condition for existence of a contractible Hamiltonian cycle in edge graph of
an equivelar triangulation of surfaces.

We moreover show that the contractible Hamiltonian cycle bounds a triangulated $2$-disk.
If the equivelar triangulation of a surface is on $n$ vertices then this disk has
exactly $n - 2$ triangles and all of its $n$ vertices lie on the boundary cycle.
We begin with some definitions.

\section{Definitions and Preliminaries}

\begin{defn} A path $P$ in a graph $G$ is a subgraph $P:[v_1, v_2, \ldots, v_n]$ of $G$,
such that the vertex set of $P$ is $V(P) = \{v_1, v_2, \ldots, v_n\}$ and
$v_{i}v_{i + 1}$ are edges in $P$ for $1 \leq i \leq n - 1$.
\end{defn}

\begin{defn} A path $P:[v_1, v_2, \ldots, v_n]$ in $G$ is said to be a cycle if $v_{n}v_1$
is also an edge in $P$.

\end{defn}

\begin{defn} A graph without any cycles or loops is called a tree

\end{defn}

If a surface $S$ has an equivelar triangulation on $n$ vertices then the proof of the
Theorem \ref{T2} is given by considering a tree with $n-2$ vertices in the dual map of
the degree-regular triangulation of the surface. We define this tree as follows\,:

\begin{defn} Let $M$ denote a  map on a surface $S$, which is the dual map of
a $n$ vertex degree-regular triangulation $K$ of the surface. Let $T$ denote a tree on
$n - 2$ vertices on $M$. We say that $T$ is a proper tree if\,:
\begin{enumerate}
\item whenever two vertices $u_1$ and $u_2$ of $T$ belong to a face $F$ in $M$, a path
$P[u_1, u_2]$ joining $u_1$ and $u_2$ in boundary of $F$ belongs to $T$.

\item any path $P$ in $T$ which lies in a face $F$ of $M$ is of length at most $q - 2$,
where $M$ is a map of type $\{q, 3\}$.
\end{enumerate}
\end{defn}

If $v$ is a vertex of a simplicial complex $X$, then the number of edges containing $v$ is
called the degree of $v$ and is denoted by $\deg_X(v)$ (or $\deg(v)$). If the number of $i$-simplices
of a simplicial complex $X$ is $f_i(X)$ ($0\leq  i\leq 2$), then the number
$\chi(X) = f_0(X) - f_1(X) + f_2(X)$ is called the {\em Euler characteristic of X}. A simplicial complex
is called neighbourly if each pair of its vertices form an edge.

\section{Example : An equivelar-triangulation and its dual}
\hrule
\medskip
\begin{center}
{\large\bf {Non-Orientable degree-regular combinatorial 2-manifold
of $\chi = -2$.}}
\end{center}
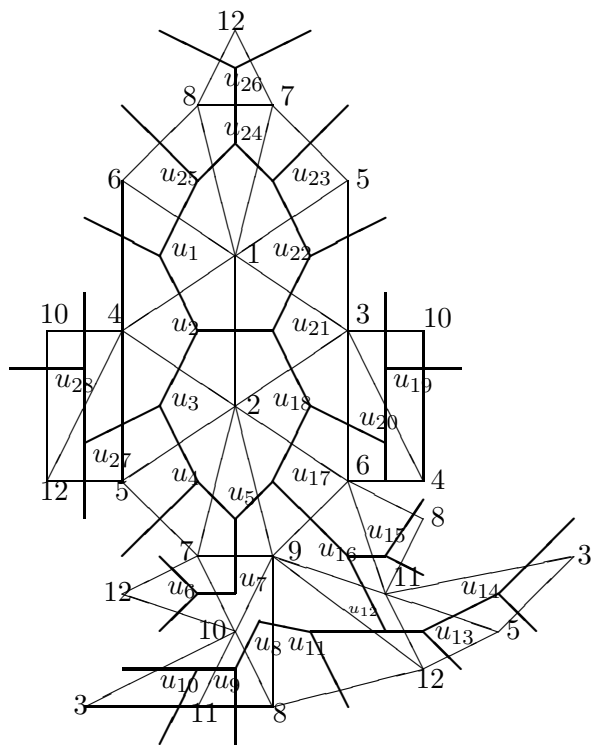
\setlength{\unitlength}{5mm}

\begin{picture}(20,28)(0,0)



\put(10,14){\line(0,1){8}}
\put(10,14){\line(1,-1){2}}
\put(12,12){\line(1,0){2}}
\put(13,16){\line(1,-4){1}}
\put(13,16){\line(-1,-4){1}}
\put(13,16){\line(-3,-2){3}}

\put(16,14){\line(0,1){8}}

\put(16,14){\line(-1,-1){2}}
\put(13,16){\line(3,-2){3}}
\put(10,22){\line(3,-2){3}}
\put(16,22){\line(-3,-2){3}}
\put(10,18){\line(3,-2){3}}
\put(13,20){\line(1,4){1}}
\put(13,20){\line(-1,4){1}}
\put(10,22){\line(1,1){2}}
\put(16,22){\line(-1,1){2}}
\put(12,24){\line(1,0){2}}

\put(13,16){\line(0,1){4}}
\put(13,16){\line(3,2){3}}
\put(13,20){\line(3,-2){3}}
\put(13,20){\line(-3,-2){3}}

\put(16,14){\line(1,-3){1}}
\put(12,12){\line(1,-2){1}}
\put(14,12){\line(-1,-2){1}}
\put(17,11){\line(-3,1){3}}
\put(13,10){\line(1,-2){1}}
\put(17,11){\line(1,-2){1}}
\put(14,8){\line(0,1){4}}
\put(18,9){\line(-4,3){4}}
\put(14,8){\line(4,1){4}}

\put(16,14){\line(2,-1){2}}
\put(18,13){\line(-1,-2){1}}
\put(12,12){\line(-2,-1){2}}
\put(13,10){\line(-3,1){3}}
\put(13,10){\line(-1,-2){1}}
\put(14,8){\line(-1,0){2}}
\put(12,8){\line(-1,0){3}}
\put(13,10){\line(-2,-1){4}}
\put(18,9){\line(2,1){2}}
\put(17,11){\line(3,-1){3}}
\put(20,10){\line(1,1){2}}
\put(17,11){\line(5,1){5}}

\put(12,24){\line(1,2){1}}
\put(14,24){\line(-1,2){1}}
\put(10,14){\line(-1,0){2}}
\put(10,18){\line(-1,0){2}}
\put(8,18){\line(0,-1){4}}
\put(10,18){\line(-1,-2){2}}

\put(16,14){\line(1,0){2}}
\put(16,18){\line(1,0){2}}
\put(18,18){\line(0,-1){4}}
\put(16,18){\line(1,-2){2}}
\put(9.8,13.5){\mbox{\small 5}} \put(7.8,13.5){\mbox{\small 12}}
\put(7.8,18.2){\mbox{\small 10}}\put(9.6,18.2){\mbox{\small 4}}
\put(9.6,21.8){\mbox{\small 6}}\put(11.6,24){\mbox{\small 8}}
\put(12.5,26){\mbox{\small 12}}\put(14.2,24){\mbox{\small 7}}
\put(16.2,21.8){\mbox{\small 5}}\put(16.2,18.2){\mbox{\small 3}}
\put(16.2,14.2){\mbox{\small 6}}\put(18.2,13.8){\mbox{\small 4}}
\put(18,18.1){\mbox{\small 10}}\put(18.2,12.8){\mbox{\small 8}}
\put(17.2,11.2){\mbox{\small 11}}
\put(13.3,19.8){\mbox{\small 1}}\put(13.3,15.8){\mbox{\small 2}}
\put(11.5,11.9){\mbox{\small 7}}\put(14.4,11.9){\mbox{\small 9}}
\put(9.5,10.8){\mbox{\small 12}}\put(12,9.8){\mbox{\small 10}}
\put(8.7,7.8){\mbox{\small 3}}\put(11.8,7.6){\mbox{\small 11}}
\put(14,7.6){\mbox{\small 8}}\put(17.8,8.5){\mbox{\small 12}}
\put(20.2,9.8){\mbox{\small 5}}\put(22.1,11.8){\mbox{\small 3}}

\put(11.3,18){\bf \small $u_2$}\put(14.5,18){\bf \small $u_{21}$}
\put(11.3,20){\bf \small $u_1$}\put(14,20){\bf \small $u_{22}$}
\put(11.3,16){\bf \small $u_3$}\put(14,16){\bf \small $u_{18}$}
\put(11,22){\bf \small $u_{25}$}\put(14.5,22){\bf \small $u_{23}$}
\put(12.7,23.2){\bf \small $u_{24}$}\put(12.7,24.5){\bf \small $u_{26}$}
\put(9.2,14.5){\bf \small $u_{27}$}\put(16.3,15.5){\bf \small $u_{20}$}
\put(8.2,16.5){\bf \small $u_{28}$}\put(17.2,16.5){\bf \small $u_{19}$}
\put(11.3,14){\bf \small $u_4$}
\put(12.8,13.5){\bf \small $u_5$}\put(13.1,11.2){\bf \small $u_7$}
\put(15.2,12){\bf \small $u_{16}$}
\put(16.5,12.5){\bf \small $u_{15}$}\put(14.5,14){\bf \small $u_{17}$}
\put(11.2,11){\bf \small $u_6$}\put(11,8.5){\bf \small $u_{10}$}
\put(12.4,8.5){\bf \small $u_9$}\put(13.5,9.5){\bf \small $u_8$}
\put(14.4,9.5){\bf \small $u_{11}$}\put(16,10.5){\bf \tiny $u_{12}$}
\put(18.3,9.8){\bf \small $u_{13}$}\put(19,11){\bf \small $u_{14}$}

\thicklines

\put(12,18){\line(1,0){2}}
\put(12,18){\line(-1,2){1}}
\put(14,18){\line(1,2){1}}
\put(11,20){\line(1,2){1}}\put(11,20){\line(-2,1){2}}
\put(15,20){\line(-1,2){1}}\put(15,20){\line(2,1){2}}
\put(12,22){\line(-1,1){2}}
\put(14,22){\line(1,1){2}}
\put(12,22){\line(1,1){1}}\put(14,22){\line(-1,1){1}}
\put(13,23){\line(0,1){2}}\put(13,25){\line(-2,1){2}}\put(13,25){\line(2,1){2}}

\put(12,18){\line(-1,-2){1}}\put(14,18){\line(1,-2){1}}
\put(11,16){\line(1,-2){1}}\put(15,16){\line(-1,-2){1}}
\put(12,14){\line(1,-1){1}}\put(14,14){\line(-1,-1){1}}
\put(13,13){\line(0,-1){2}}
\put(11,16){\line(-2,-1){2}}\put(15,16){\line(2,-1){2}}
\put(9,15){\line(0,1){2}}\put(17,15){\line(0,1){2}}
\put(9,15){\line(0,-1){2}}\put(17,15){\line(0,-1){1}}
\put(9,17){\line(0,1){2}}\put(9,17){\line(-1,0){2}}
\put(17,17){\line(0,1){2}}\put(17,17){\line(1,0){2}}
\put(12,14){\line(-1,-1){2}} \put(14,14){\line(1,-1){2}}
\put(13,11){\line(-1,0){1}}\put(13,11){\line(1,-1){0.65}}
\put(13,9){\line(1,2){0.65}}\put(13,9){\line(-1,0){1}}
\put(13,9){\line(0,-1){2}}\put(12,9){\line(-1,0){2}}
\put(12,9){\line(-1,-2){1}}

\put(15,10){\line(-5,1){1.35}}\put(15,10){\line(1,0){2}}
\put(15,10){\line(1,-2){1}}
\put(17,10){\line(1,0){1}}\put(18,10){\line(1,-1){1}}
\put(18,10){\line(2,1){2}}
\put(20,11){\line(1,-1){1}}\put(20,11){\line(1,1){2}}
\put(12,11){\line(-1,-1){1}}\put(12,11){\line(-1,1){1}}

\put(16,12){\line(1,0){1}}
\put(16,12){\line(1,-2){1}}\put(17,12){\line(2,3){1}}
\put(17,12){\line(2,-1){1}}

\end{picture}


\hrule
\smallskip

\section{Some facts about proper trees.}

\begin{lemma} Let $v \in V(T)$ be a vertex in a proper tree $T$. Then $\deg(v) \leq 3$.

\end{lemma}

\noindent \begin{proof} Let $M$ denote the dual map of a triangulation $K$. Thus
$\deg(u) = 3$ for all $u \in V(M)$. Since $T$ is a subgraph of the edge graph of
$M$, $\deg(v) \leq 3$ for all $v \in V(T)$.
\end{proof}\hfill$\Box$

\smallskip

\begin{lemma} Let $T$ be a proper tree and $m$ be the number of vertices of degree 3 in $T$.
Then the number of vertices of degree one in $T$ = $m + 2$.

\end{lemma}

\smallskip

\noindent \begin{proof} Let $P_1$ denote a path of maximum length in $T$. Then $P_1$ has two
ends which are also ends of $T$, for otherwise $P_1$ will not be of maximum length.
If there is no vertex of degree 3 in $T$ which also lies in $P_1$ then
$P_1 = T$, as $T$ is connected, and we are done. Otherwise, let $u_1$ be a vertex of degree 3
such that $u_1 \in P_1 \cap T$. Let $u_1$ be the initial point of a path $P_2$ of maximum length in the tree
$T' = T \setminus P_1$. Thus $P_2$ is edge disjoint with $P_1$. Then by the above argument the  end of $P_2$ other than
$u_1$ is also an end of $T'$ and hence of $T$. Further, if there is a vertex $w_1$ of degree 3 on $P_2 \cap T'$, we repeat
the above process to find an end of $T$. Thus,
for each vertex of degree 3 in $T$ we get an end of $T$. This together with ends of $P_1$ proves
that the number of ends of $T =  m + 2$.
\end{proof}\hfill$\Box$

\smallskip

\begin{lemma} Let $T$ be a proper tree in a polyhedral map $M$ of type $\{q, 3\}$
on a surface $S$. Then $T\bigcap F \neq \emptyset$ for any face $F$ of $M$.
\end{lemma}

\smallskip

\noindent \proof\,\, Let $e$ denote the number of vertices of degree one in $T$. Since $T$ has $n - 2$ vertices,
it has $n - 3$ edges. We claim that the $n - 3$ edges of $T$ lie in exactly $n - e$ faces of $M$.

To prove this we enumerate the number of faces of $M$ with which the edges of $T$ are incident.

We construct sets $E$ and $\tilde{F}$ as follows. Let $E$ be a singleton set which contains an edge
$e_1$ of $T$ and $F_1$ and $F_2$ be the faces of $M$ such that $e_1$ lies in them. Put
$\tilde{F} := \{F_1, F_2\}$.
Add an adjacent edge $e_2$ of $e_1$ to $E$. There is exactly one face $F_3$ different from $F_1$ and
$F_2$ such that $e_2$ lies in $F_3$. Add this to set $\tilde{F}$ to obtain $\tilde{F} := \{F_1, F_2, F_3\}$.
Successively, we add edges to the set $E$ which are adjacent to edges in $E$ till we
exhaust all the edges of $T$. Each additional edge added to $E$ contributes exactly one face
to the set $\tilde{F}$ unless it is adjacent to two edges in the set $E$. Thus the number of
faces in $\tilde{F} =$ (number of edges of $T$ - number of vertices of degree three) + 1. In
a 3-tree, the number of vertices of degree 3 = number of end point - 2. Thus $\#\tilde{F} =
n - 3 - (e - 2) + 1$. That is $\#\tilde{F} = n - e$.

Let $F(M)$ denote the set of all faces of $M$. Let $G = F(M) \setminus \tilde{F}$. Then $\# G = e$.
We claim that an end vertex of $T$ lies on exactly one face $F \in G$. Observe that each vertex
$u$ of $T$ is incident with exactly three distinct faces $F_1$, $F_2$ and $F_3$ of $M$. The edge
of $T$ incident with $u$ lies in two of these faces, say $F_1$ and $F_2$, $i.e.$, $F_1$, $F_2
\in \tilde{F}$. Since, $u$ is an end vertex, there is no edge of $T$ which is incident with
$F_3$, for otherwise this violates the definition of $T$. Thus $u$ is incident with exactly
one face $F_3$ of $M$ such that $F_3 \in G$. Since, $u$ is an arbitrary end point this hypothesis
holds for all the end vertices. If it happens that for some end vertices $u_1$ and $u_2$ of $T$,
the corresponding faces $W_1 = W_2 \in G$ then we would have $u_1$ and $u_2$ on the same face
of $M$ but no path on $W_1$ joining $u_1$ and $u_2$ lies in $T$. This contradicts the definition
of $T$. Thus $G$ has exactly $e$ distinct elements. This proves the lemma. \hfill $\Box$


\smallskip

\begin{lemma}\label{l4}$\!\!$ Let $K$ be a $n$ vertex degree regular triangulation of a surface $S$. Let $M$
denote the dual polyhedron corresponding to $K$ and $T$ be a $n - 2$ vertex proper tree in
$M$. Let $D$ denote the subcomplex of $K$ which is dual of $T$. Then $D$ is a triangulated
2-disk and $bd(D)$ is a Hamiltonian cycle in $K$.
\end{lemma}

\noindent\proof\, By definition of a dual, $D$ consists of $n - 2$ triangles corresponding to
$n - 2$ vertices of $T$. Two triangles in $D$ have an edge in common if the corresponding
vertices are adjacent in $T$. It is easy to see that $D$ is a collapsible simplicial complex
and hence it is a triangulated 2-disk.

Moreover, since $T$ has vertices of degree one, $bd(D) \neq \emptyset$, and being boundary
complex of a 2-disk it is a connected cycle. Observe that
the number of edges in $n - 2$ triangles is $3(n - 2)$
and for each edge of $T$ exactly $2$ edges are identified. Hence the number of edges which
remain unidentified in $D$ is $3(n - 2) - 2(n - 3) = n$. Similarly the number of vertices in
$bd(D)\colon = \partial{D} = n$. If there are vertices $v_1, v_2 \in \partial{D}$ such that $v_1$ and $v_2$ lie
on a path of length $< n$ and $v_1 = v_2$. This means there are faces $F_1$ and $F_2$ in $D$
with $v_1 \in F_1$, $v_2 \in F_2$, $F_1 \neq F_2$ and $F_1$ not adjacent to $F_2$. Thus there
exist a face $F'$ in $D$ such that the vertex $u_{F'}$ in $T$ corresponding to $F'$, 
does not belong to the face $F(v_1)$ corresponding to vertex $v_1$. But
this contradicts that $T$ is a proper tree. Thus all the cycle $\partial{D}$ contains
exactly $n$ distinct vertices. Since $\#V(K) = n$, $\partial{D}$ is a Hamiltonian cycle in $K$. \hfill $\Box$

\begin{theo} \label{T2}The edge graph $EG(K)$ of an equivelar triangulation $K$ of a surface has a contractible
Hamiltonian cycle if and only if the edge graph of corresponding dual map $M$ of $K$ has a proper tree.
\end{theo}

\begin{proof} The above Lemma \ref{l4} shows the if part.
Conversely, let $K$ denote an equivelar triangulation and $H := (v_1, v_2, v_3, \ldots, v_n)$ denote
a contractible Hamiltonian cycle in $EG(K)$. Let $\tau_1, \tau_2, \ldots, \tau_m$ denote the faces
of triangulated disk whose boundary is $H$. We claim that all the triangles have their vertices on
boundary of the disk, i.e. on $H$. For otherwise there will be identifications on the surface because
all the vertices of $K$ also lie on $H$. If $x$ denotes the number of triangles in this disk then the
Euler characteristic relation gives us $1 = n - [\frac{(3\times x) - n}{2} + n] + x$. Thus, $x = n - 2$.
So that $m = n - 2$. Now, in the edge graph of dual map $M$ of $K$, consider the graph corresponding to this
disk whose vertices correspond to the dual of faces $\tau_1, \tau_2, \ldots, \tau_m$.
Now it is easy to check that this graph is a tree which is also a proper tree.
\end{proof}\hfill $\Box$
\bigskip

\noindent {\bf  Acknowledgement\,:} The author thanks D. Barnette \cite{barnette2} for reading
and appreciating the idea of Proper Tree. That this tree may be a necessary and
sufficient condition for existence of separating Hamiltonian cycle (Theorem \ref{T2}) was
suggested by him.


{\small

\begin{thebibliography}{99}
\bibitem{a1}
A. Altshuler, Construction and enumeration of regular maps
on the torus, {\em Discrete Math.} ({\bf 4}) (1973), 201--217.
\bibitem{a2}
A. Altshuler, Hamiltonian circuits in some maps on the torus,
{\em Discrete Math.} ({\bf 4}) vol. 1, (1972), 299--314.
\bibitem{barnette1}
Barnette, D.: 3-Trees in Polyhedral Maps, {\em Israel J. Math.}
{\bf 79}, (1992), 251 - 256.
\bibitem{barnette2}
Barnette, D.: Personal Communications
\bibitem{barnette3}
Barnette, D.: Contractible circuits in 3-connected graphs, {\em Discrete Math.}
{\bf 187}, (1998), 19 - 29.
\bibitem{BondyMurthy}
J. A. Bondy and U. S. R. Murthy, {\em Graph theory with applications},
North Holland, Amsterdam, 1982.
\bibitem{dn}
B. Datta and N. Nilakantan, Equivelar polyhedra with few vertices,
{\em Discrete \& Comput Geom.} {\bf 26} (2001), 429--461.
\bibitem{duke}
Duke, R. A.: On the Genus and Connectivity of Hamiltonian Graphs,
{\em Discrete Math.}, {\bf 2}, (1972), 199 - 206.
\bibitem{grunbaum}
Gr$\ddot{u}$nbaum, B.: Polytopes, graphs and complexes, {\em Bull. Amer. Math. Soc.}, {\bf 76}, (1970), 1131 - 1201.
\bibitem{m}
J. R. Munkres, {\em Elements of Algebraic Topology\/},
Addison-Wesley, California, 1984.
\bibitem{PulpakaVince1}
Pulpaka, H. and Vince, A.: Non-revisiting Paths on Surfaces with Low Genus,
{\em Discrete Math.}, {\bf 182}, (1998), 267 - 277
\bibitem{PulpakaVince2}
Pulpaka, H. and Vince, A.: Non-revisiting Paths on Surfaces,
{\em Discrete Comput. Geom.}, {\bf 15}, (1996), 353 - 357
\bibitem{Tutte}
Tutte, W. T.: A theorem on planar graphs, {\em Trans. Amer. Math. Soc.}, {\bf 82}, (1956), 99 - 116.
\bibitem{yu}
X. Yu, Disjoint paths, planarizing cycles and spanning walks,
{\em Trans. Amer. Math. Soc.} ({\bf 4}) vol 349, (1997), 1333--1358.
\end{thebibliography}
}
\end{document}